\documentclass[12pt]{amsart}
\usepackage{graphicx,amsfonts}
\usepackage{amsmath}
\usepackage[all]{xy}
\newcommand{\g}{\mathfrak{g}}
\newcommand{\C}{\mathbb{C}}
\newcommand{\R}{\mathbb{R}}
\newcommand{\F}{\mathbb{F}}
\newcommand{\Q}{\mathbb{Q}}
\newcommand{\gl}{\mathfrak{gl}}
\newcommand{\A}{\mathcal{A}}
\newcommand{\B}{\mathcal{B}}
\newcommand{\Hom}{\mathop{\rm Hom}\nolimits}
\renewcommand{\phi}{\varphi}

\newtheorem*{definition}{Definition}
\newtheorem{lemma}{Lemma}
\newtheorem{proposition}{Proposition}
\newtheorem*{theorem}{Theorem}

\title{Detecting the orientation of long links by finite type invariants}
\author{S.V.Duzhin, M.V.Karev}

\begin{document}
\begin{abstract}
We prove the existence of a degree 7 Vassiliev invariant of long (or string)
two-component links which is not preserved under the simultaneous change of
orientation of both components. The non-invertibility of this invariant
can be detected by the standard weight system with values in the tensor
square of the universal enveloping algebra for $\gl(n)$.
\end{abstract}

\maketitle

\section{Introduction}

It is well known that classical knot polynomials and quantum invariants in
general take equal values on a knot and its inverse. The class of Vassiliev
invariants is strictly wider \cite{Vog}, and the problem whether they can be
used to tell a knot from its inverse is open and seemingly difficult. In
this paper, we discuss the corresponding problem for links with more than 1
component.

Let $S^1_p$ be the disjoint union of $p$ numbered copies of an oriented
circle, $\R^1_p$ the disjoint union of $p$ numbered copies of an oriented
real line, and $I^1_p$ the disjoint union of $p$ numbered copies of an
oriented interval $[0,1]$.

\begin{definition}
An $p$-component closed link is a smooth embedding of $S^1_p$ into oriented
3-space $\R^3$, considered up to a component-preserving
isotopy.

An $p$-component long link is a smooth embedding of $R^1_p$ into oriented
3-space $\R^3$ with fixed behavior at infinity: $x_i(t)=[i,0,t]$ for
$|t|>C$, considered up to an
isotopy, fixed at infinity, and a Euclidean movement of space.

An $p$-component string link is a smooth embedding of $I^1_p$ into the
slice $\R^2\times[0,1]\subset\R^3$ with fixed endpoints: $x_i(t)=[i,0,t]$ 
for $t=0,1$, considered up to an
isotopy, fixed at the boundary, and a Euclidean movement of the slice.
\end{definition}

\begin{center}
\includegraphics[height=35mm]{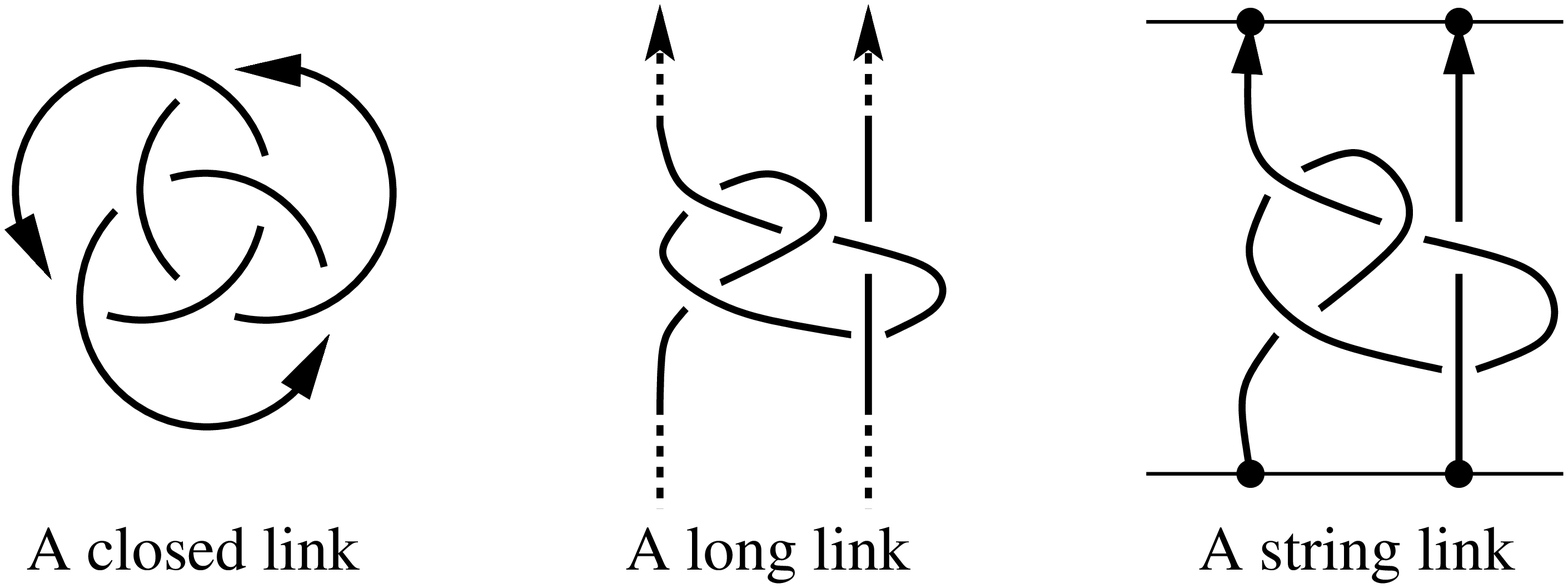}
\end{center}

It is immediately clear that the theories of long and string links are
equivalent. In the case of knots ($p=1$) they are also equivalent to the
theory of closed links; this is not true when $p>1$.

The \textit{problem of invertibility} is stated as follows.
Given a link $L$, let $L'$ be the same link with the orientation of all
components reversed. May it happen that $L'$ is not equivalent to $L$?
If yes, then what invariants can be used to detect that difference?

For knots, the problem was standing for a long time, until Trotter
\cite{Tro} proved in 1964 that some knots (e.g. pretzel knot with
parameters $(3,5,7)$) are non-invertible. The simplest non-invertible knot is
$8_{17}$ (\cite{Ka}). However,
it is not known if the two knots $8_{17}$ and $8_{17}'$, or any other pair
of mutually inverse knots, can be distinguished
by finite type invariants. 

In the case of links, the only published result is a theorem of X.-S.\,Lin
\cite{Lin2} that asserts the non-invertibility of Vassiliev invariants for
closed links with at least 6 components. Other papers containing the words
`links', `Vassiliev invariants', `invertibility' (`separation', `detection'
etc) in their titles/abstracts have smaller relevance to the plain problem
we are speaking about, e.g. D.Bar-Natan \cite{BN2} studies \textit{homotopy}
string links invariants, X.-S.Lin \cite{Lin1} uses a modified definition of
finite type invariants, etc.

The present text deals with the problem of detecting the orientation by
finite type invariants for \textit{long links}. The reformulation of the
problem in terms of chord diagrams shows immediately that the answer is
positive for $p>2$. In the case of 2-component links ($p=2$) the problem is
nontrivial; below, we give a proof of the following

\begin{theorem}
There exists a Vassiliev invariant $f$ of order no greater than 7
and a 2-component long link $L$ such that $f(L')\not= f(L)$.
\end{theorem}

We will give two proofs of this theorem, both based on direct calculations.
The first proof (Proposition \ref{prop1}, page \pageref{prop1}) uses chord
diagrams and requires computer calculations. The second one (Proposition
\ref{prop2}, page \pageref{prop2}) uses Jacobi diagrams and can be done by
hand. Frankly speaking, both of these proofs refer to the invariants of
\textit{framed} links. In Section \ref{defram} we explain why they also
imply the same result in the unframed case.

%A reader with enough background does not have to read the proofs --- it suffices to look
%at the non-invertible diagrams on pp.~\pageref{prop1} and \pageref{prop2},
%respectively.

\section{Reduction to chord diagrams}

Finite type invariants for various types of links are defined in the same
way as they are in the classical case of closed knots, see \cite{BN1,Sta}. 

Let $\F$ be a field of characteristic 0, e.g. $\Q$ or $\C$. Denote by
$V_n(p)$ the vector space of $\F$-valued Vassiliev invariants for
$p$-component long links of degree no greater than $n$, by $\A_n(p)$ the
space spanned by chord diagrams of degree $n$ on $p$ strings modulo 4-term
relations, and by $W_n(p)=\A_n^*(p)=\Hom_\Q(\A_n(p),\F)$ the corresponding
space of \textit{weight systems}.
Then there is a linear map
$\sigma_n^p:V_n(p)\to W_n(p)$ (taking the \textit{symbol} of a Vassiliev
invariant) whose kernel coincides with $V_{n-1}(p)$ and whose image
consists of all weight systems that vanish on chord diagrams with
isolated chords. For \textit{framed} links, the image equals the entire 
$W_n(p)$, and, until Section \ref{defram}, we will actually be speaking about
invariants of framed links.

As before, we denote by $L'$ the inverse of a link $L$. Setting 
$\tau_V(f)(L)=f(L')$ for $f\in V$, we obtain an involution $\tau_V$ 
in the space of Vassiliev invariants. On the level of chord diagrams, 
the corresponding operation $\tau_A$ acts as follows:
$$
\tau_A :\quad 
\raisebox{-9mm}{\includegraphics[height=20mm]{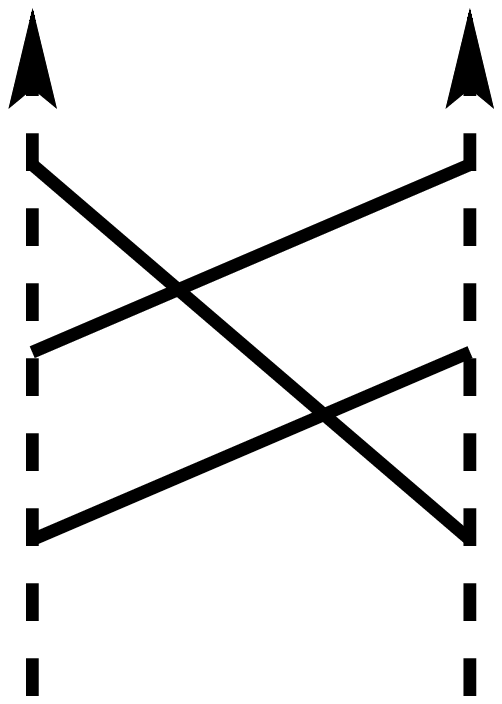}}
\quad\longmapsto\quad
\raisebox{-9mm}{\includegraphics[height=20mm]{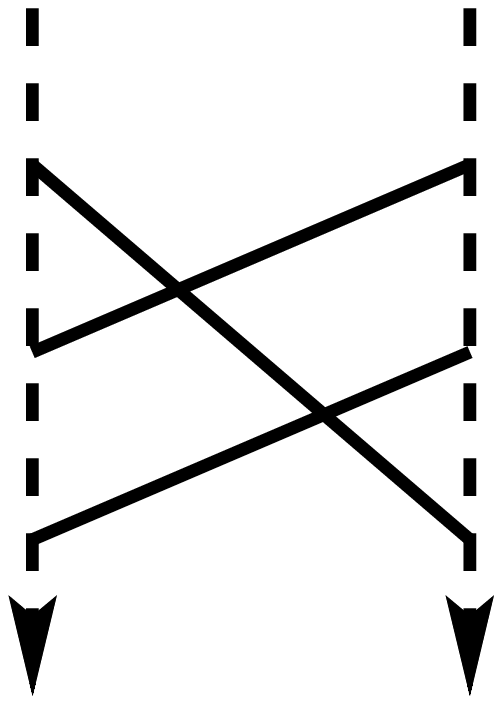}}
\quad=\quad
\raisebox{-9mm}{\includegraphics[height=20mm]{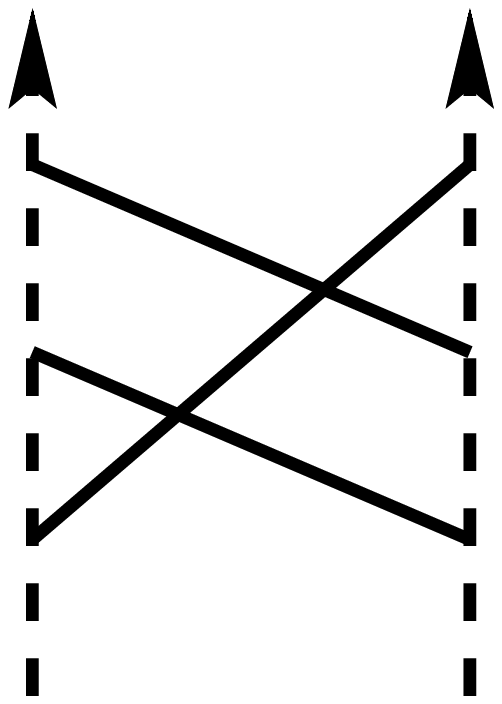}}\ ,
\label{cd3}
$$
i.e. it changes the orientation of all the strings of a chord
diagram, or, which is the same thing, reflects the plane picture
of the diagram in a horizontal line, assuming that the support is drawn
vertically (and preserves its orientation).

Instead of usual chord diagrams only, the same space $\A(p)$ may be 
spanned by all generalized chord diagrams (\textit{Chinese character
diagrams} in the terminology of \cite{BN1}), i.e. 1-3-valent graphs,
attached to the same support ($\R^1_p$) with a
cyclic order of edges specified in every 3-valent vertex.

Generalized chord diagrams
can be understood as linear combinations of ordinary chord diagrams
by iterative applications of the STU rule, e.g.
$$
\raisebox{-8mm}{\includegraphics[height=17mm]{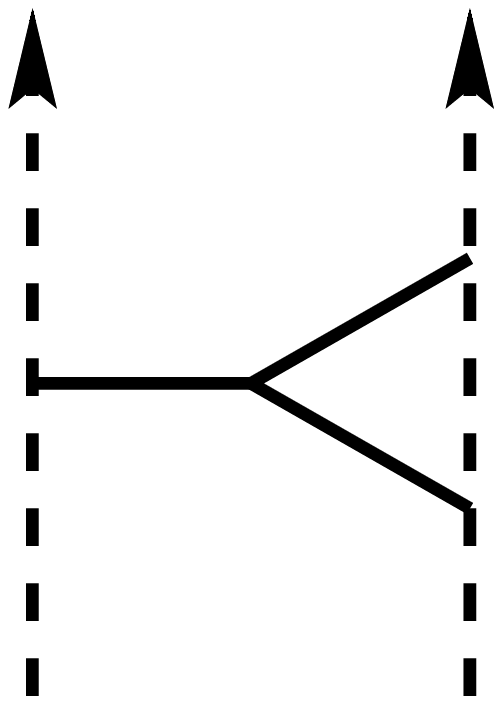}}
\quad\longmapsto\quad
\raisebox{-8mm}{\includegraphics[height=17mm]{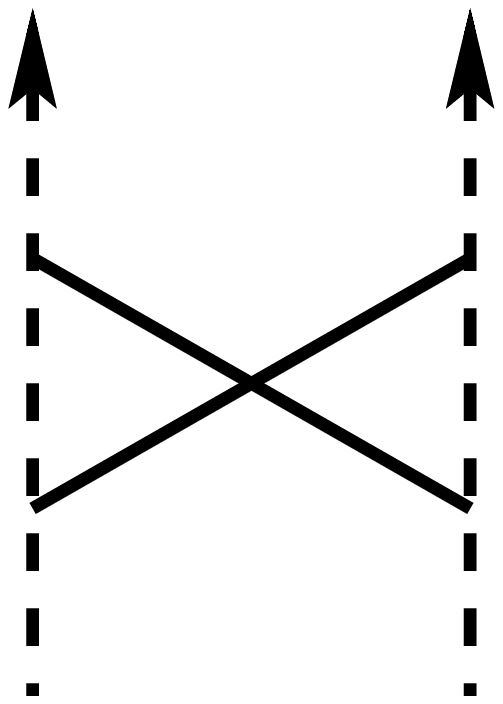}}
\quad-\quad
\raisebox{-8mm}{\includegraphics[height=17mm]{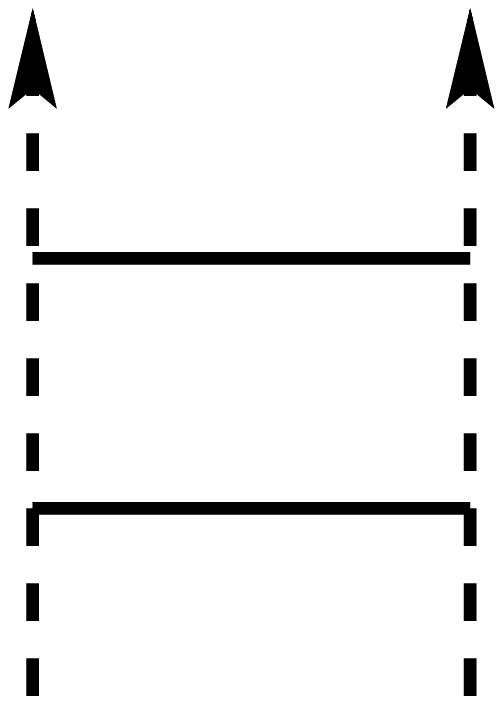}}
$$
(here and below, we use the \textit{blackboard convention}: edges around
trivalent vertices for every diagram drawn on paper are always ordered
counterclockwise).

For generalized diagrams, one should be more cautious with the definition
of the involution $\tau_A$. In fact, extending the previously defined
map $\tau_A$ from chord diagrams to generalized diagrams by STU relations,
one arrives at the following rule:
$$
\tau_A :\quad 
\raisebox{-12mm}{\includegraphics[height=25mm]{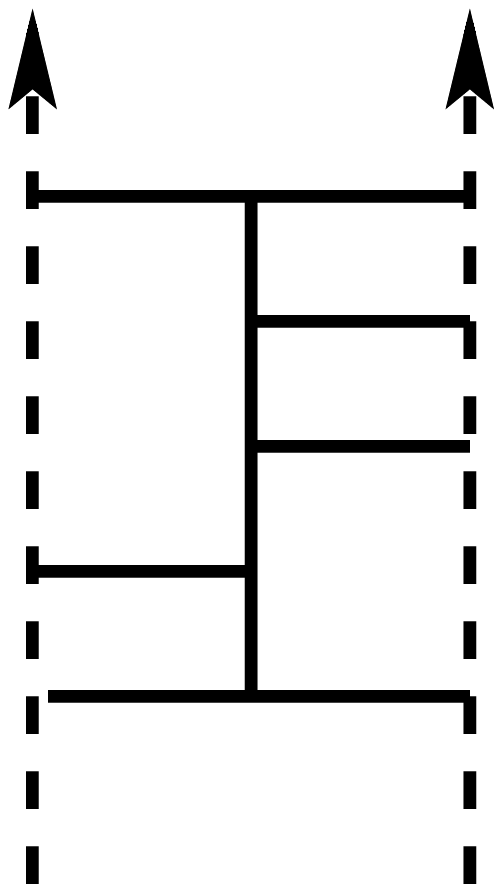}}
\quad\longmapsto\quad
-\quad\raisebox{-12mm}{\includegraphics[height=25mm]{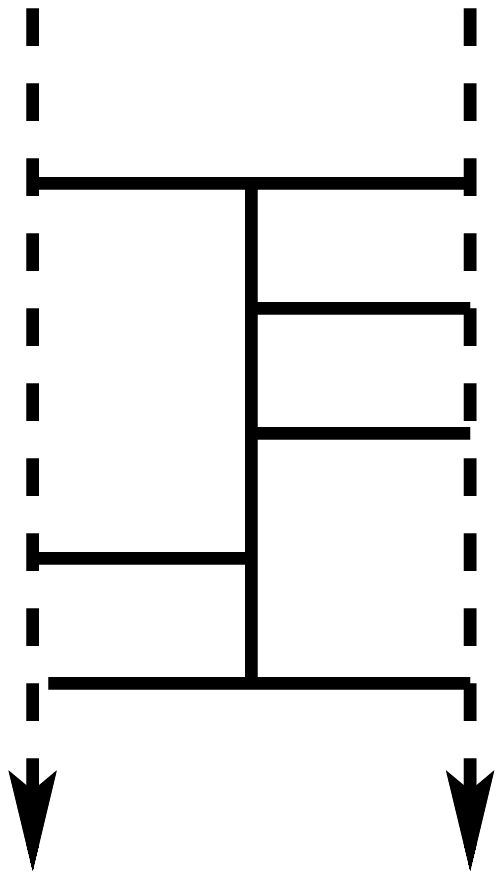}}
\quad=\quad
\raisebox{-12mm}{\includegraphics[height=25mm]{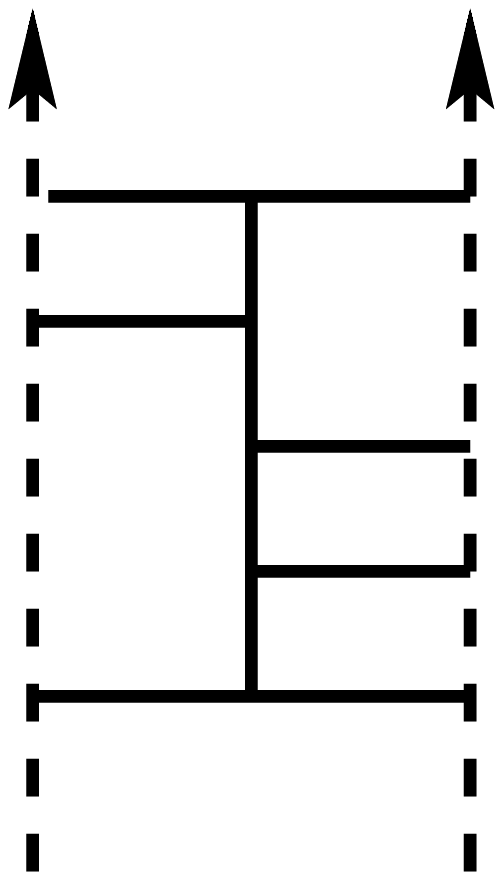}}
$$
i.e. one must either change the orientation of all strings \textit{and
multiply the result by $\pm1$} depending on the parity of the number of
trivalent vertices, or, again, simply reflect the plane picture
of the diagram in a horizontal line.

\begin{lemma}
The involution $\tau_A$ is the graded part of the involution $\tau_V$, i.e.
the following square commutes:
$$
\xymatrix{
V_n(p) \ar[r]^{\sigma_n^p} \ar[d]^{\tau_V} & W_n(p) \ar[d]^{\tau_A^*} \\
V_n(p) \ar[r]^{\sigma_n^p} & W_n(p) 
}
$$
\end{lemma}

\begin{proof}
This follows immediately from the definition
of $\sigma$ (\cite{BN1,CD}).
\end{proof}

The problem of detecting the inversion of orientation on long links
by finite type invariants is thus equivalent to the problem: \textit{find a
chord diagram on $p$ lines which is different from its inverse modulo
4T relations}. If $p\ge3$, then there is a straightforward example:
$$
\tau_A :\quad 
\raisebox{-7mm}{\includegraphics[height=17mm]{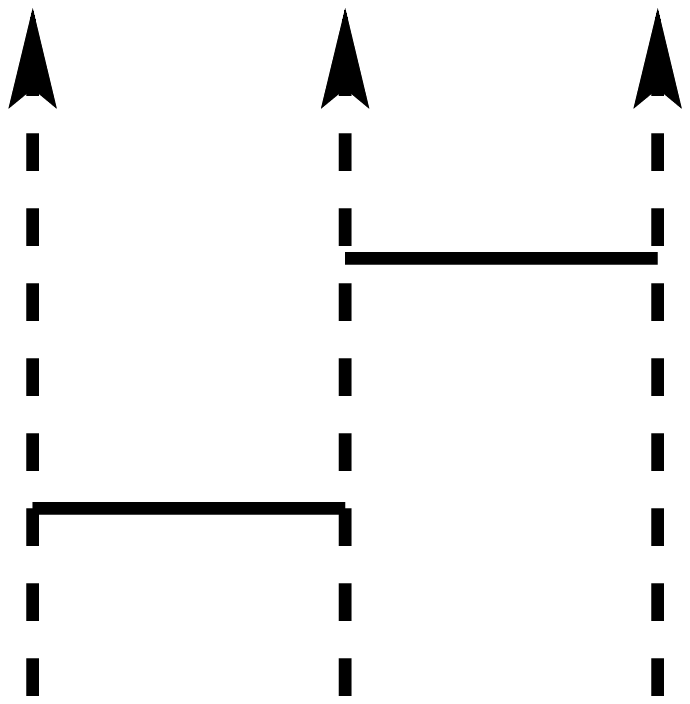}}
\quad\longmapsto\quad
\raisebox{-7mm}{\includegraphics[height=17mm]{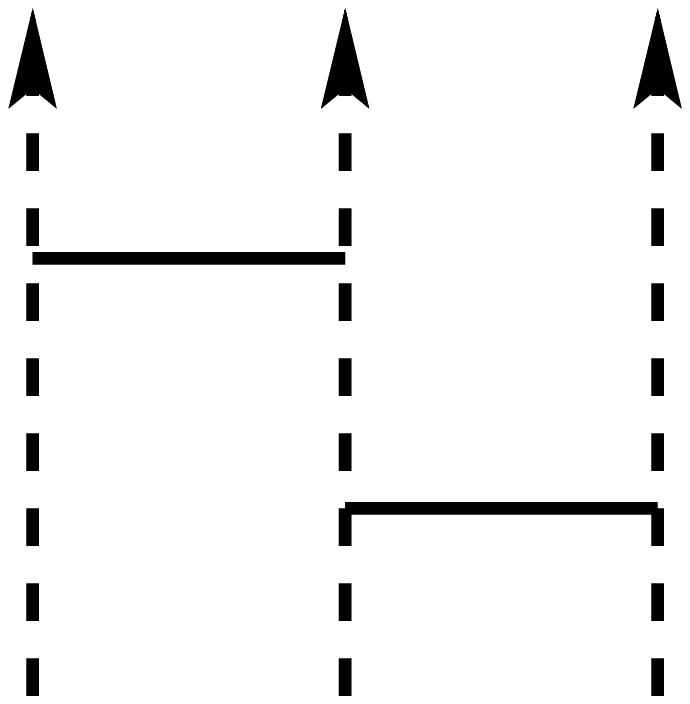}}
\quad\not=\quad
\raisebox{-7mm}{\includegraphics[height=17mm]{AB.eps}}
$$

For $n=2$ the problem is non-trivial, because, as one can check,
chord diagrams of small degrees (like the one depicted on page
\pageref{cd3}) are all $\tau_A$-invariant.
In the next section, we indicate a diagram which is not. 

\section{First proof of the theorem}
\label{secphi}

To prove that a certain element of the space $\A(p)$ is non-zero, one can use 
\textit{weight systems}, i.e. linear functionals on these spaces.
A powerful weight system is provided by Kontsevich's 
homomorphism $\phi=\phi_\g:\A(p)\to U(\g)^{\otimes p}$ for a metrized 
Lie algebra $\g$ (see \cite{Kon,CD}). In fact, $\phi$ takes values
in the $\g$-invariant subalgebra $U(p)=[U(\g)^{\otimes p}]^\g$.
We give an explicit description for the case
$\g=\gl_N$, using the basis $e_{ij}$ (matrix with one 1 and many 0's)
and the metric defined by the conjugacy rule $e_{ij}^*=e_{ji}$.

\begin{lemma}
Let $D\in\A(p)$ be a (generalized) chord diagram on $p$ strings. The element
$\phi(D)\in U(\gl_N)^{\otimes p}$ can be obtained as follows. Consider the
alternating sum of all resolutions \cite{BN1} of the inner triple points of 
the diagram $D$. For each resolution label the connected components of the 
obtained diagram by different independent indices; then replace each pair of
adjacent indices by $e_{ij}$ and take the sum over all indices from 1 to
$N$. (Closed components, if any, turn into factors $N$.)
\end{lemma}

The proof repeats the proof for the particular case $p=1$ that can be found
in \cite{BN1,CD}. We give an example instead:
$$
\raisebox{-10mm}{\includegraphics[height=21mm]{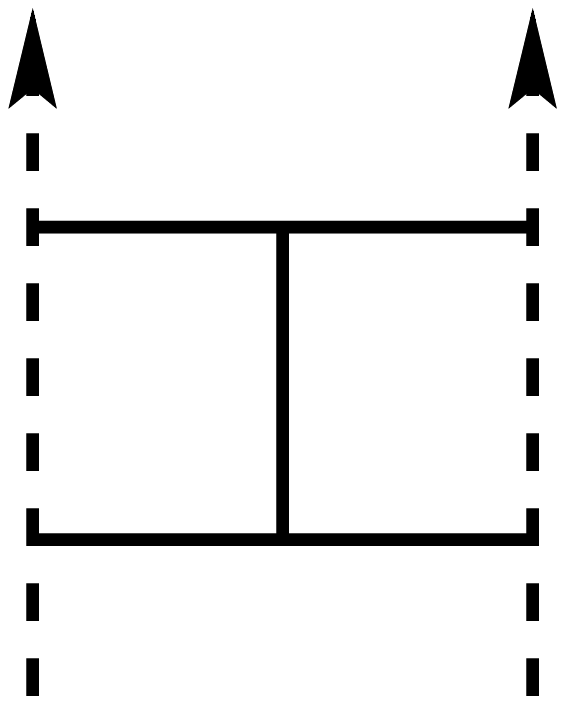}}
\quad\longmapsto\quad
\raisebox{-10mm}{\includegraphics[height=21mm]{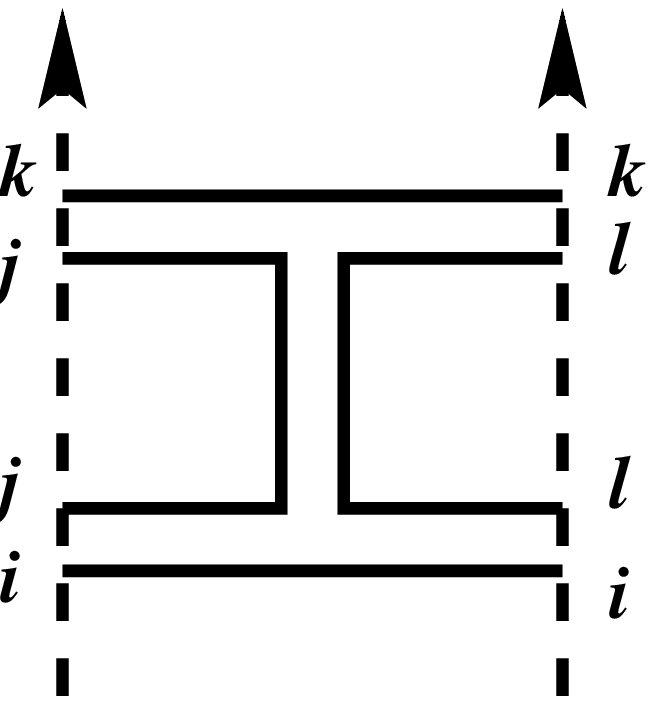}}
-
\raisebox{-10mm}{\includegraphics[height=21mm]{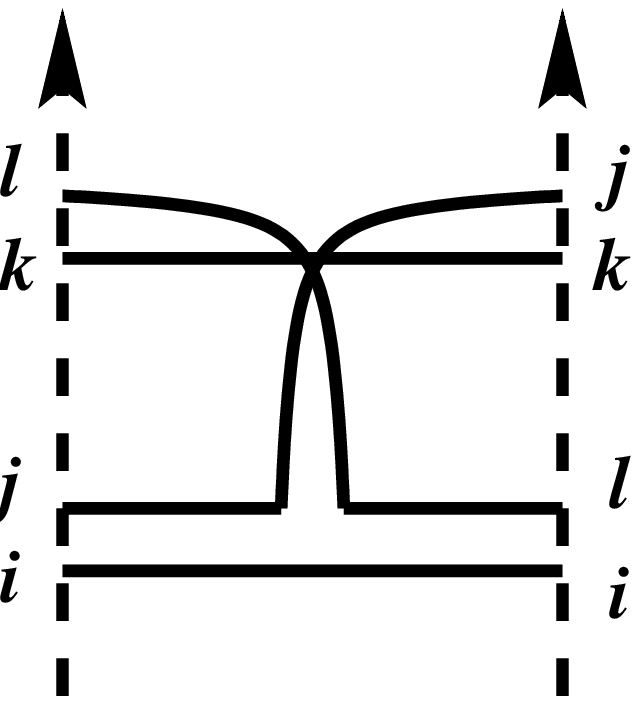}}
-
\raisebox{-10mm}{\includegraphics[height=21mm]{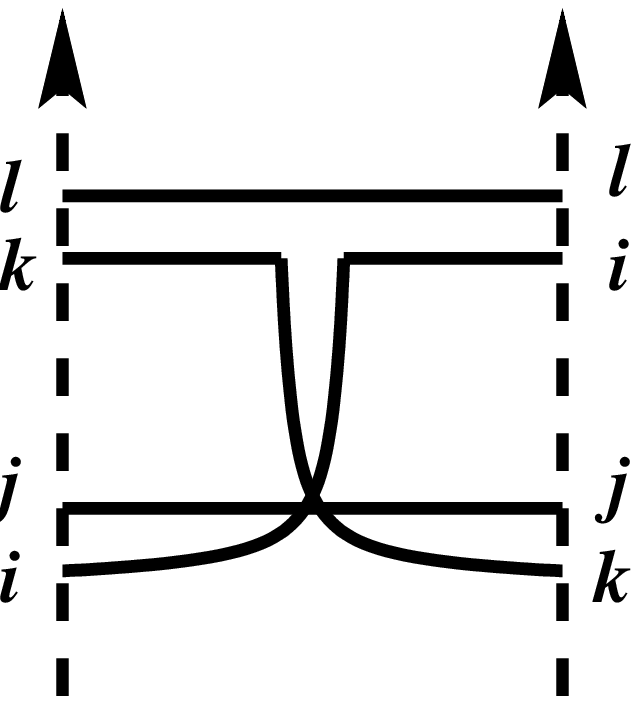}}
+
\raisebox{-10mm}{\includegraphics[height=21mm]{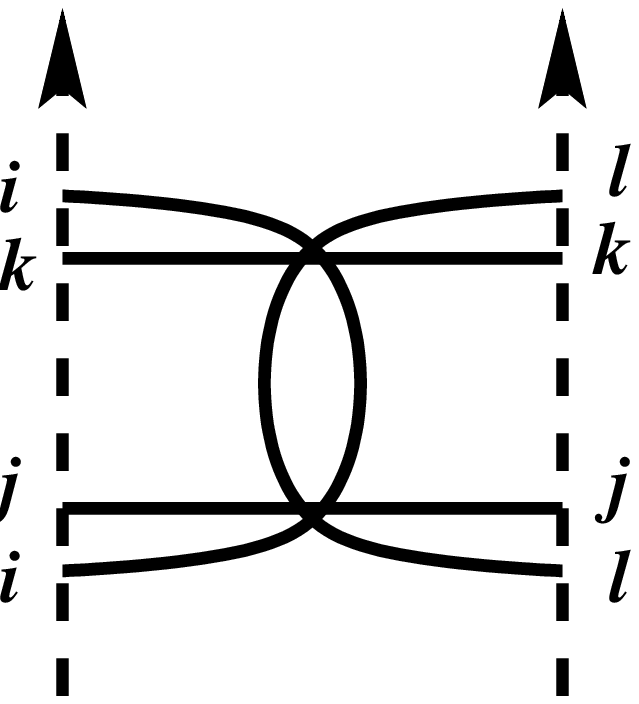}}
$$
Therefore, the image of this diagram under $\phi$ is
$$
\sum_{i,j,k,l=1}^N(e_{ij}e_{jk}\otimes e_{li}e_{kl}
-e_{ij}e_{kl}\otimes e_{li}e_{jk}
-e_{ij}e_{kl}\otimes e_{jk}e_{li}
+e_{ij}e_{ki}\otimes e_{jl}e_{lk}).
$$
Note that the order of the factors $e_{ij}$ agrees with the orientation
of each component of the support, but the order of the two subscripts
$(i,j)$ in each factor follows a go-round pattern, in our case bottom-to-top
for the left line and top-to-bottom for the right line. 

Denote by $\tau_U$ the operator in $U(\gl_N)^{\otimes p}$ that
rewrites each monomial backwards and changes the order of the two
subscripts in each generator $e_{ij}$, e.g.
$\tau_U(e_{12}e_{23}\otimes e_{13}e_{24})
=e_{32}e_{21}\otimes e_{42}e_{31}$. This is an involution that preserves
the ad-invariant subspace $U(p)=[U(\gl_N)^{\otimes p}]^\g$.

The construction of $\phi$ shows that the following assertion holds.

\begin{lemma}
The square
$$
\xymatrix{
\A(p) \ar[r]^{\phi} \ar[d]^{\tau_A} & U(p) \ar[d]^{\tau_U} \\
\A(p) \ar[r]^{\phi} & U(p) 
}
$$
is commutative.
\end{lemma}

It follows that the non-invertibility of a chord diagram can be checked
on the level of the universal enveloping algebra: if the $\phi$-image 
of a diagram is not $\tau_U$-invariant, then the diagram itself is not
invariant under $\tau_A$.
Therefore, the next proposition gives the first proof of our
Theorem.

\begin{proposition}
Each of the following diagrams is different from its image under the
involution $\tau_A$:
\begin{center}
\raisebox{-12mm}{\includegraphics[height=25mm]{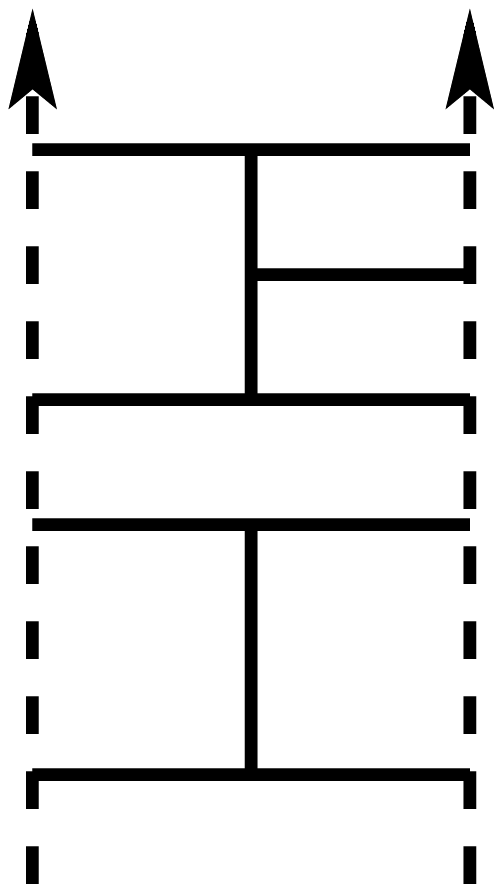}}
\quad,\qquad
\raisebox{-12mm}{\includegraphics[height=25mm]{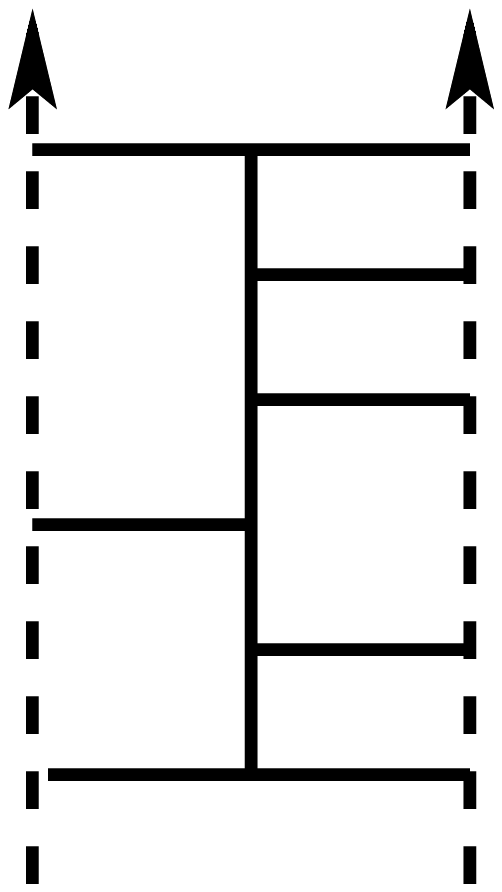}}
\quad.
\end{center}
\label{prop1}
\end{proposition}

\begin{proof}
A direct computer calculation of the $\phi$-images of each diagram and its
inverse shows that they are different for the Lie algebra $\gl_4$
(and hence for any $\gl_N$ with $n\ge4$).
The calculations are organized as follows: we fix a lexicographical
ordering of the basic elements $e_{ij}$ and transform the expression
for $\phi(D)$ obtained by the above algorithm using the 
commutator relations between the generators 
in order to rewrite each monomial lexicographically.
The programs, input and output files are available at
\cite{Calc}. The result for the left-hand side diagram, for example,
is a polynomial consisting of 58378 terms; its calculation requires
several hours on a reasonably fast modern PC.
\end{proof}

\section{Reduction to Jacobi diagrams}

Another, yet simpler, reformulation of the invertibility problem,
can be given in terms of colored \textit{Jacobi} diagrams.
A colored Jacobi diagram is the same thing as a 
\textit{Chinese character} as defined in \cite{BN1}, only its univalent
vertices are labeled by $p$ colors. The \textit{space} of colored Jacobi
diagrams $\B(p)$ is defined as a vector space formally generated by all
$p$-colored Jacobi diagrams modulo two kinds of relations: antisymmetry
and IHX.

Similarly to Theorem 8 of \cite{BN1}, one can prove that there is a
symmetrization map $\chi:\B(p)\to\A(p)$ which provides a linear
isomorphism between the two spaces (see also \cite{LM}).

The map $\chi$ is defined as follows --- we explain that in the case
$p=2$. Let $D$ be a Jacobi diagram with $k$ `legs' of color 1 and $l$ `legs'
of color 2, then $\chi(D)$ is the average of all the $k!l!$ ways to attach
1-colored legs to the first line of the support and 2-colored legs to the
second line of the support, e.g.
\begin{center}
\includegraphics[height=15mm]{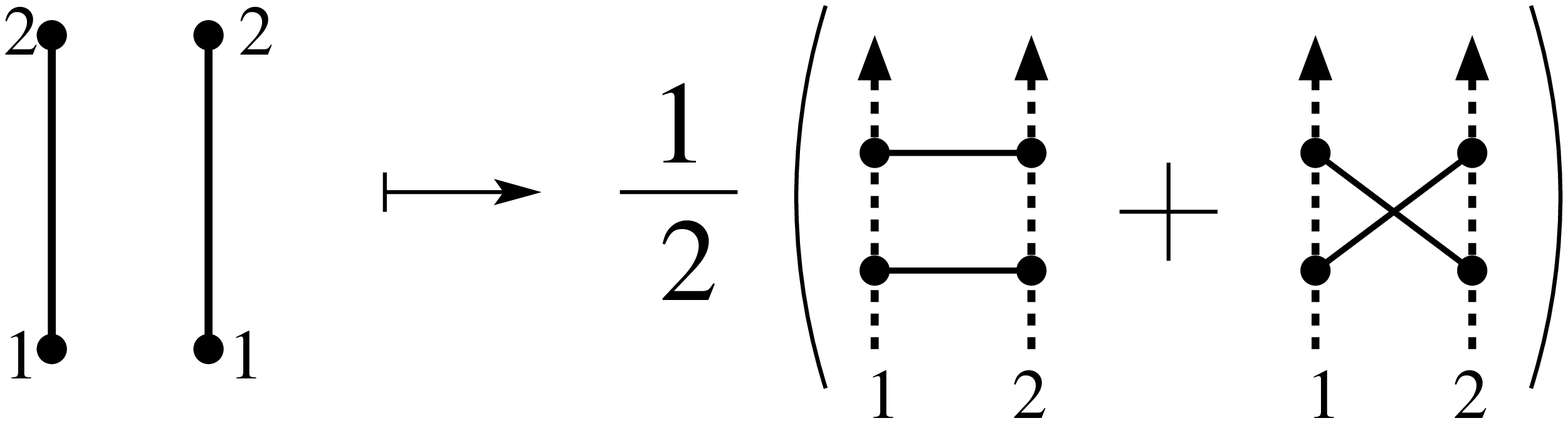}
\end{center}

Colored Jacobi diagrams can thus be viewed as symmetric elements
of the space $\A(p)$, much in the same way as commutative polynomials
over a Lie algebra $\g$ can be viewed as symmetric elements of the universal
enveloping algebra $U(\g)$, by virtue of the PBW theorem.

The isomorphism $\chi$ is very important for our needs because
the involution $\tau_A$ takes an especially simple form when transferred
to $\B(p)$ via this isomorphism. In fact, the following assertion holds:

\begin{lemma}
Let $\tau_B:\B(p)\to\B(p)$ be the linear operator defined as identity on
each Jacobi diagram with an even number of legs and as multiplication by
$-1$ on each Jacobi diagram with an odd number of legs. Then the square
$$
\xymatrix{
\B(p) \ar[r]^{\chi} \ar[d]^{\tau_B} & \A(p) \ar[d]^{\tau_A} \\
\B(p) \ar[r]^{\chi} & \A(p) 
}
$$
is commutative.
\end{lemma}

\begin{proof}
The proof smoothly follows from the definitions of both $\tau_A$ and
$\tau_B$ given above. A simple illustration might be helpful
to understand how it goes. Let
$$
  D = \raisebox{-7mm}{\includegraphics[height=15mm]{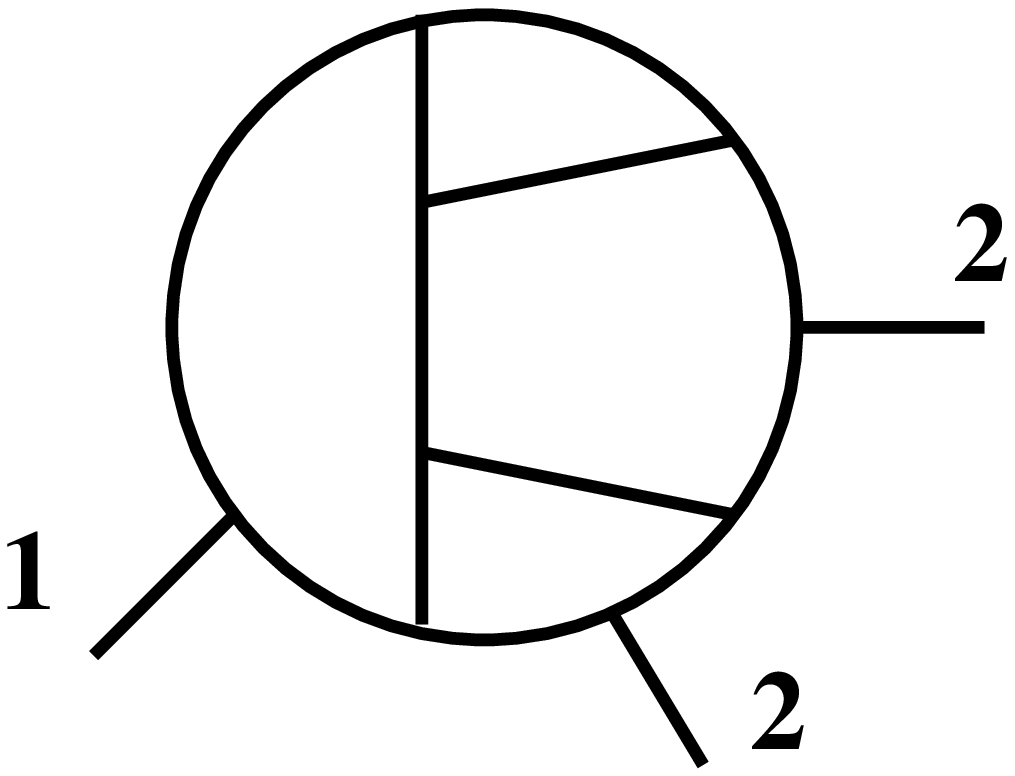}}.
$$
Then by definition
$$
  \tau_B(D) = - \raisebox{-7mm}{\includegraphics[height=15mm]{DB.eps}}.
$$
Now
$$
  \chi(D)\ =\ \frac12\ \raisebox{-7mm}{\includegraphics[height=15mm]{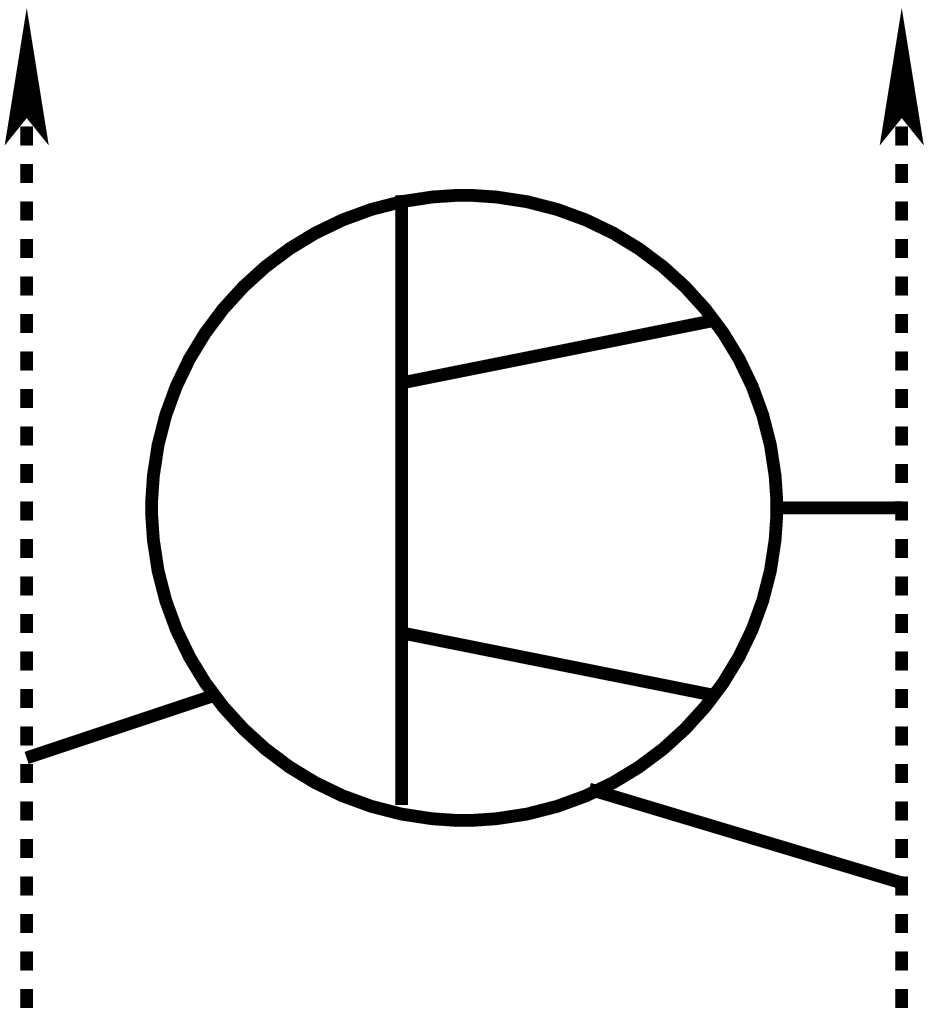}}
          \ +\ \frac12\raisebox{-7mm}{\includegraphics[height=15mm]{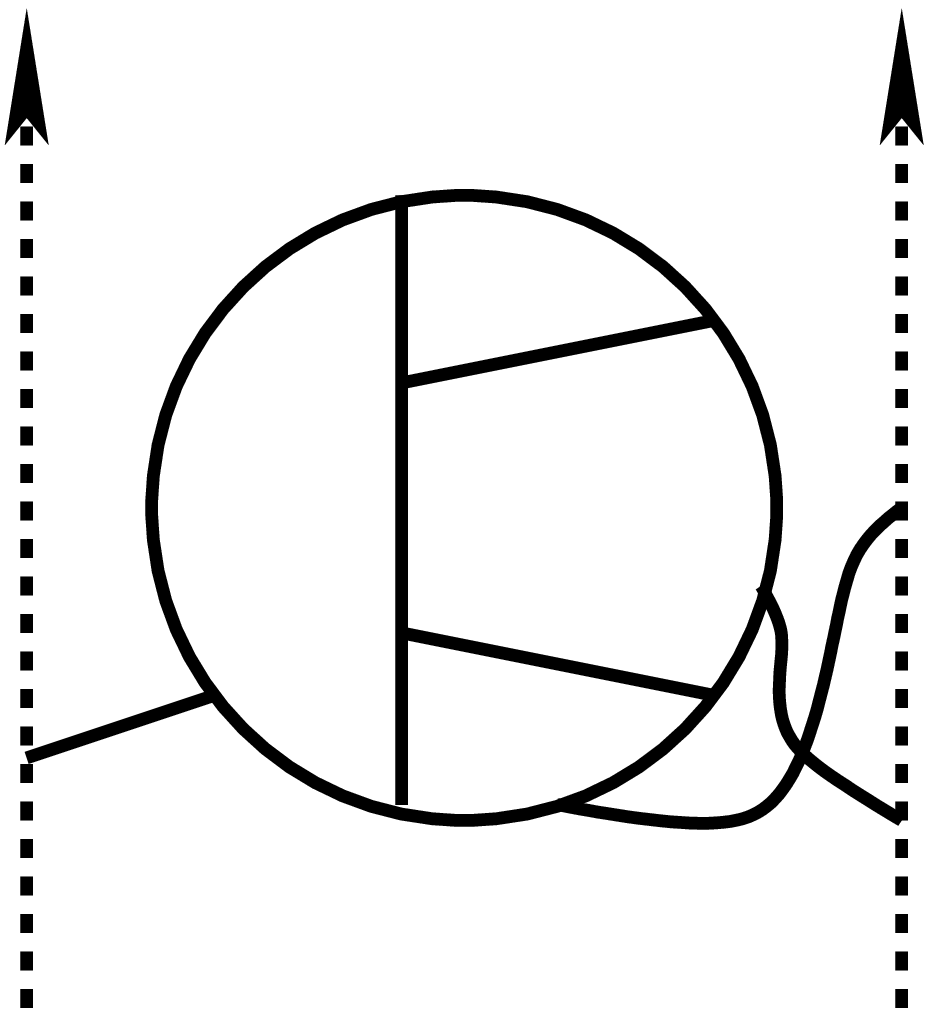}}
$$
and hence
$$
  \tau_A(\chi(D))\ =\ 
        -\frac12\raisebox{-7mm}{\includegraphics[height=15mm]{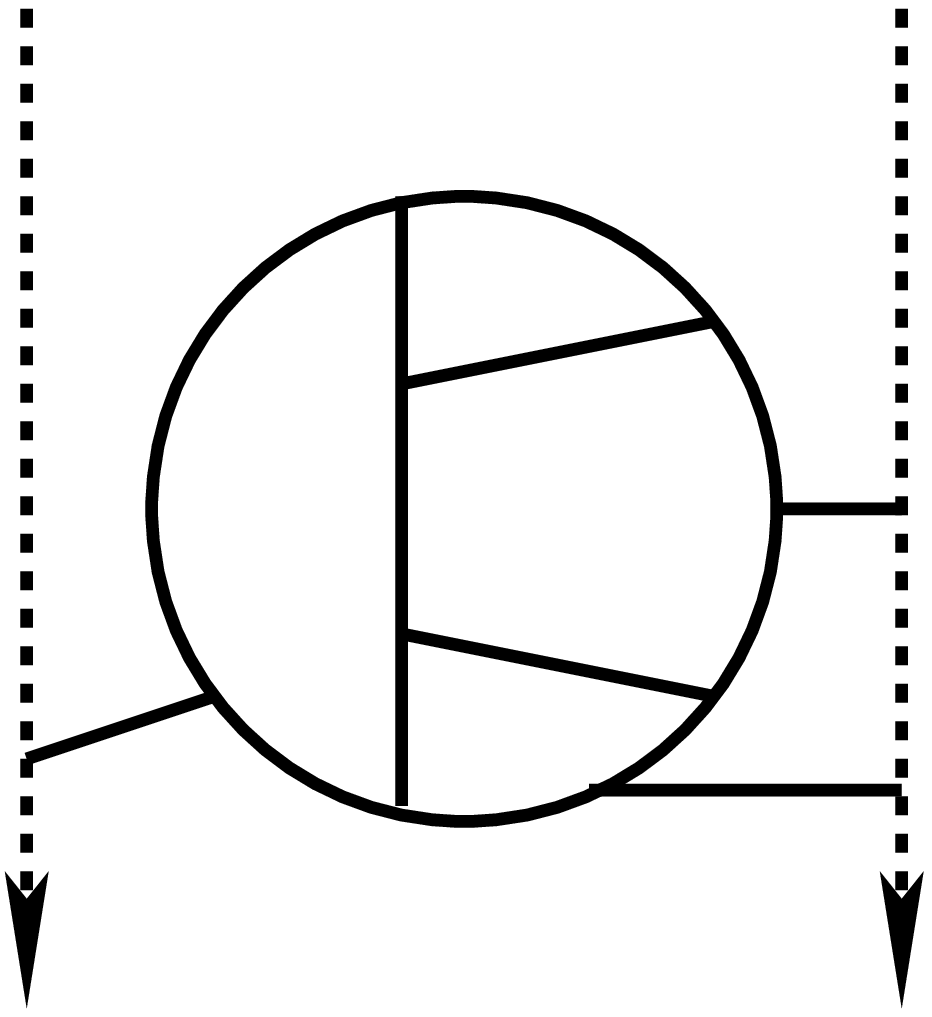}}
        -\frac12\raisebox{-7mm}{\includegraphics[height=15mm]{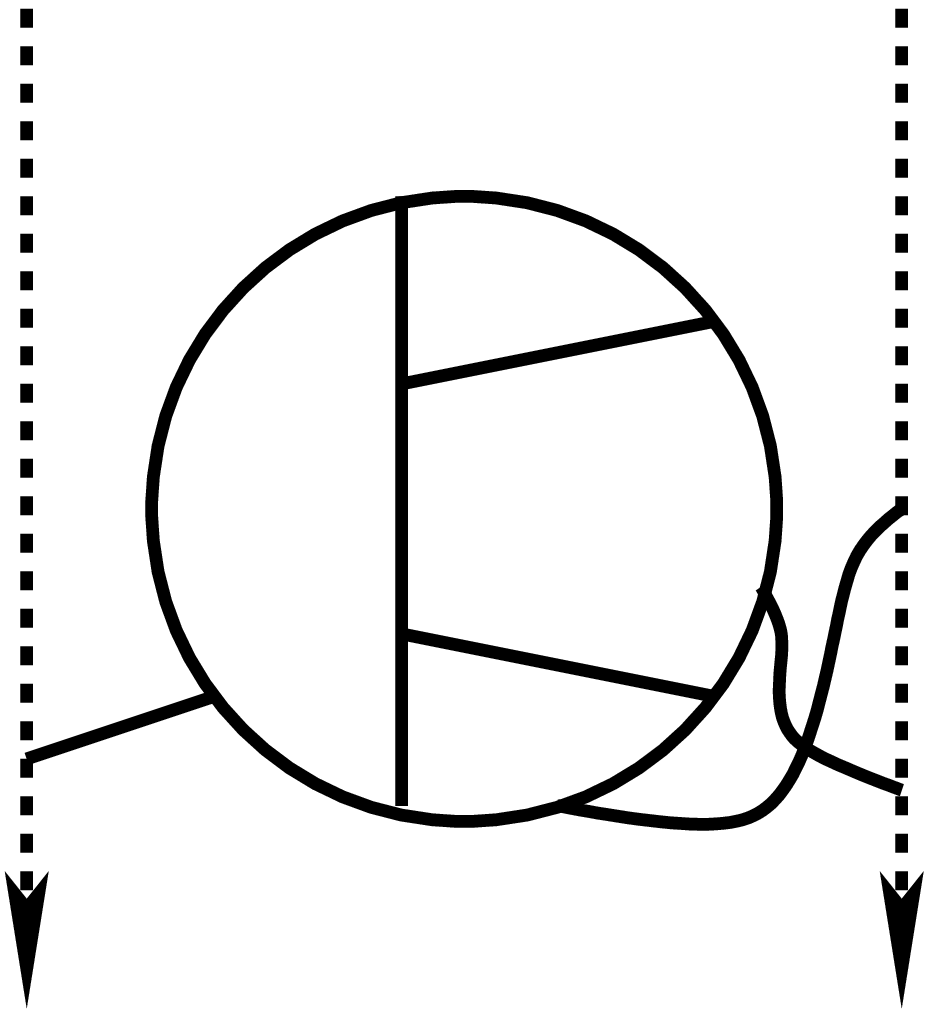}}
      \ =\ -\frac12\raisebox{-7mm}{\includegraphics[height=15mm]{DA2.eps}}
        -\frac12\raisebox{-7mm}{\includegraphics[height=15mm]{DA1.eps}}
     \ =\ \chi(\tau_B(D)).
$$
\end{proof}

The problem of finite type invertibility for 2-component long links can now
be restated as follows: \textit{is there a non-zero 2-colored Jacobi diagram
with an odd number of legs?}. In the next section, we give an example of
such a diagram.

\section{Second proof of the theorem}
\label{secpsi}

Take a metrized Lie algebra $\g$ and denote by $S(\g)$ the symmetric algebra
of the vector space $\g$. The weight system $\B(1)\to S(\g)$ described in
\cite{BN1,CD} has an immediate generalization to an arbitrary value of $p$
giving a homomorphism
$\psi:\B(p)\to S(\g)^{\otimes p}$ whose image lies in the $\g$-invariant
subalgebra $S(p)=[S(\g)^{\otimes p}]^\g$.

\begin{lemma}
The mapping $\psi$ fits into a commutative diagram
$$
\xymatrix{
\B(p) \ar[r]^{\psi} \ar[d]^{\chi} & S(p) \ar[d]^{\pi} \\
\A(p) \ar[r]^{\phi} & U(p) 
}\label{csqr1}
$$
where $\chi$ is the isomorphism defined in the previous section,
$\phi$ is the Kontsevich weight system for the algebra $\A(p)$,
and $\pi$ is the Poincar\'e-Birkhoff-Witt isomorphism raised to the
$n$-th tensor power.
\end{lemma}

\begin{proof}
The proof is straightforward, just like in the conventional case $p=1$.
\end{proof}

For the Lie algebra $\g=\gl_N$ there is a pictorial algorithm
for the calculation of $\psi$, similar to the procedure for $\phi$
described in section \ref{secphi}. The only difference is that now
the basic elements $e_{ij}$ are commutative, so we do not have to bother
about the order of factors in the monomials. Here is a simple example:
\begin{align*}
  \psi\ :\quad &\raisebox{-7mm}{\includegraphics[height=15mm]{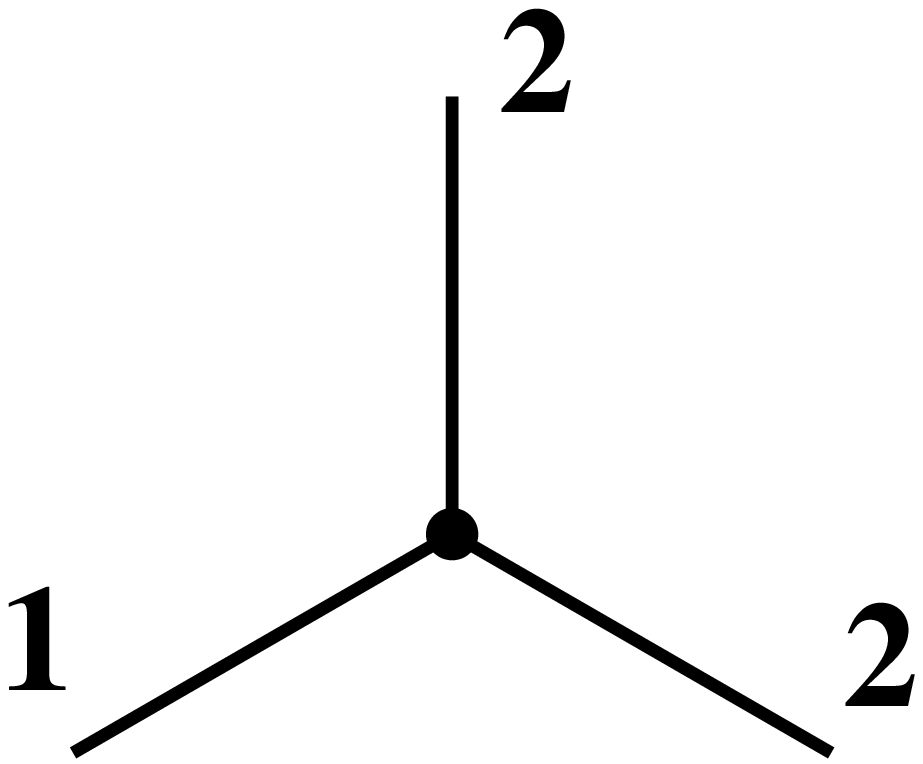}}
\quad\longmapsto\quad
       \raisebox{-9mm}{\includegraphics[height=19mm]{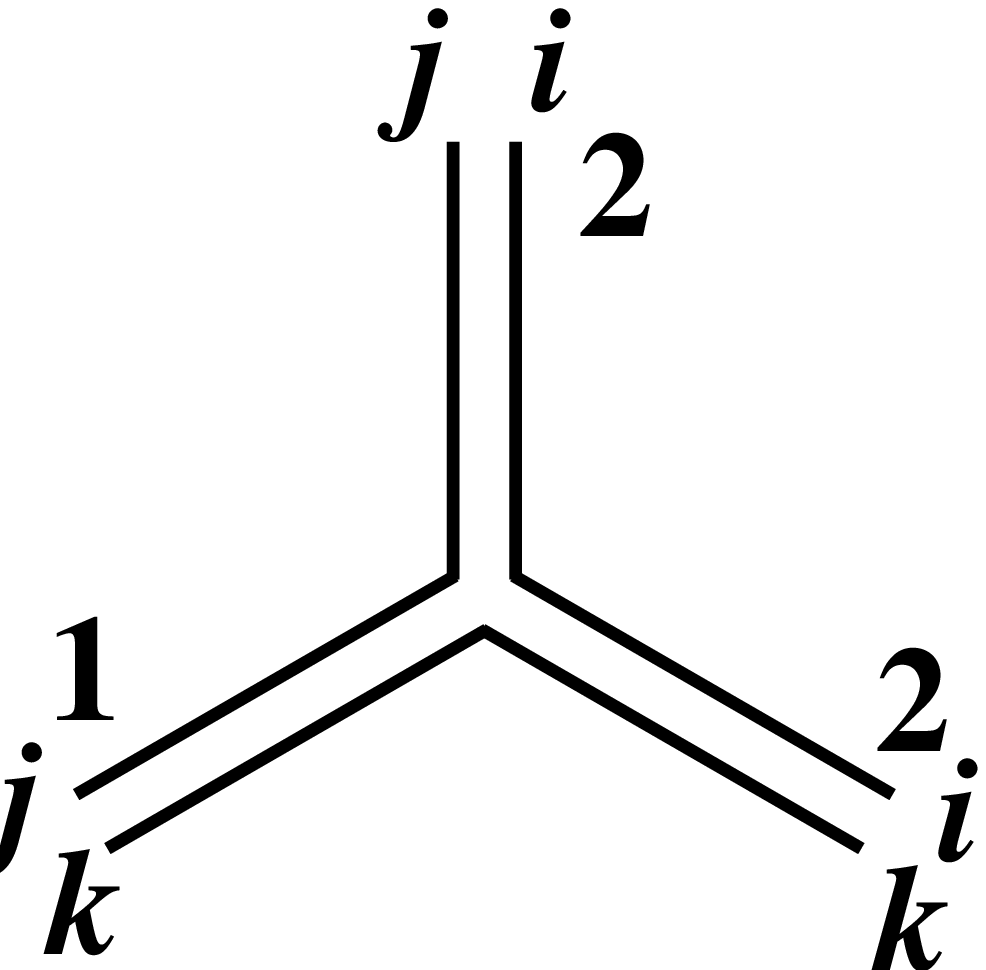}}
\quad-\quad
       \raisebox{-9mm}{\includegraphics[height=19mm]{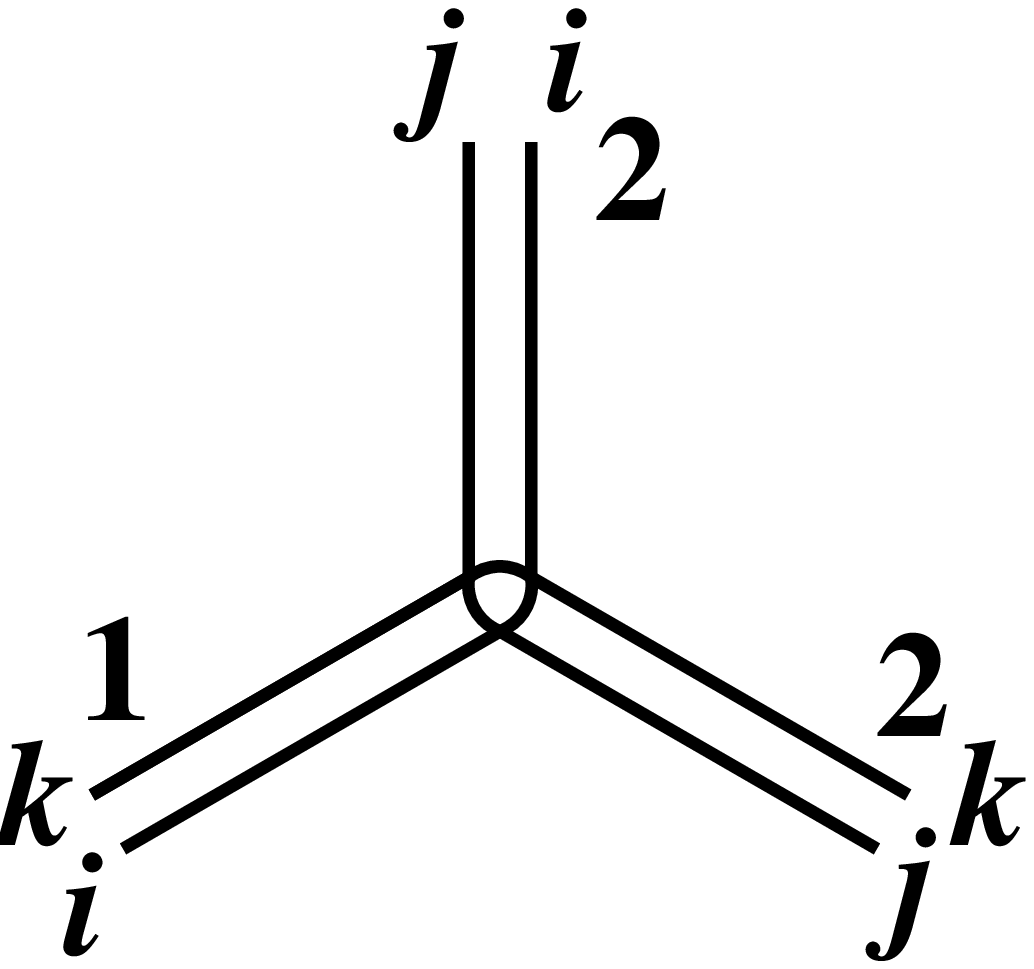}}\\
\quad&\longmapsto\quad
\sum_{i,j,k=1}^N(e_{jk}\otimes e_{ij}e_{ki}-e_{ki}\otimes e_{ij}e_{jk})
\ =\ 0.
\end{align*}

A big advantage that the mapping $\psi$ has over $\phi$ is that the 
$\g$-invariant part of $S(\gl_N)^{\otimes p}$ has a rather transparent structure.
A $p$-colored \textit{necklace} of order $n$ is a combinatorial object
defined as a sequence of $n$ numbers between 1 and $p$, considered
up to cyclic shifts; it can be best viewed as a circular arrangement
of $n$ $p$-colored beads. To every necklace one can assign an element
of $S(\gl_N)^{\otimes p}$, which turns out to be ad-invariant,
as follows. Sprinkle different indices over the arcs of the necklace;
every bead transforms into $e_{ij}$ where $i$ is the index written on the
incoming arc and $j$ is the index of the outgoing arc;
put this element $e_{ij}$ into the tensor factor whose number is the color of
the bead and take the sum over all indices from 1 to $N$. For example:
$$
[1212]\ =\ 
       \raisebox{-7mm}{\includegraphics[height=15mm]{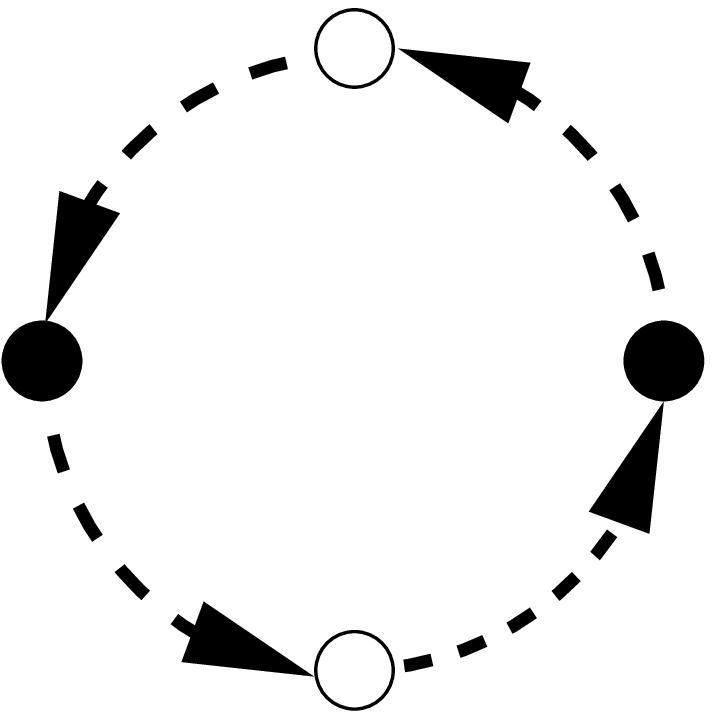}}
\quad\longmapsto\quad
       \raisebox{-7mm}{\includegraphics[height=15mm]{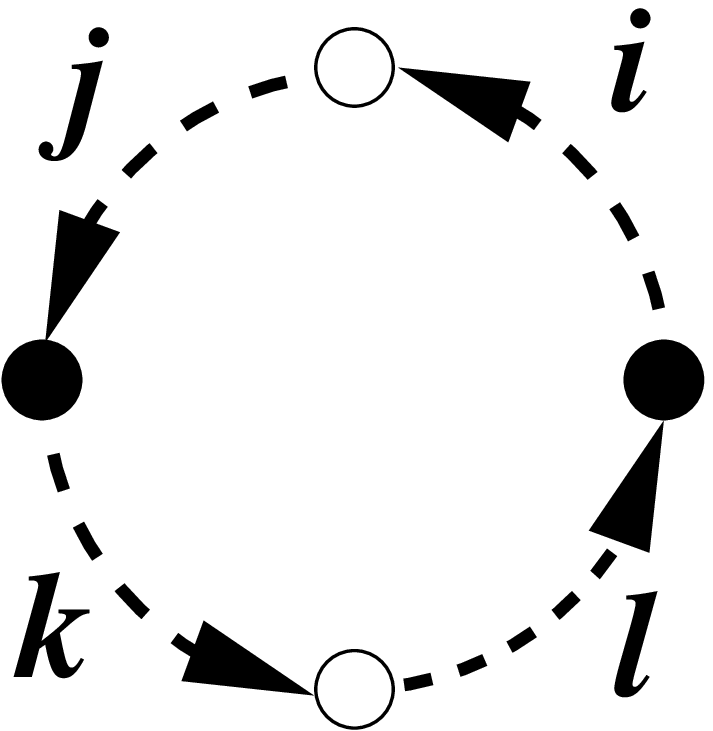}}
\quad\longmapsto\quad
\sum_{i,j,k,l=1}^N(e_{ij}e_{kl}\otimes e_{jk}e_{li})\ =:\ x_{1212}\ .
$$
We denote such \textit{necklace elements} by $x_\mu$ where $\mu$ is the
lexicographically smallest sequence corresponding to the given necklace.

\begin{lemma}
For the Lie algebra $\g=\gl_N$ the $\g$-invariant part of 
$S(\g)^{\otimes p}$, as an algebra, is generated by the necklace elements.
Algebraic relations between necklace elements of a given degree may exist
for small values of $N$, but they disappear as $N\to\infty$.
\end{lemma}

\begin{proof}
Statements close to this lemma are proved in \cite{Weyl} and \cite{VP}. 
We are sure that the fact itself must be written up somewhere,
and we would be grateful to anybody who can give us an exact reference.
\end{proof}

%Let us call the \textit{necklace algebra} the free commutative algebra
%generated by all necklaces and the variable $N$, i.e. the inductive limit
%of $[S(\gl_N)^{\otimes p}]^\g$ as $N\to\infty$.

In the algebra generated by necklaces there is an involution $\tau_S$ that
inverts the orientation of each necklace. For $p=2$, in small degrees, up
to 5, it is identical; the minimal necklace which is not fixed by $\tau_S$,
is $x_{112122}$. It turns out that the operation $\tau_S$ agrees with all
other inversions denoted by $\tau$ with various subscripts in this paper.

\begin{lemma}
The involutions of the change of orientation in the spaces
$\A(p)$, $\B(p)$, $S(p)$, $U(p)$ commute with the four arrows 
of the diagram in Lemma (\ref{csqr1}).
More precisely, we have a commutative cube
$$\xymatrix{
& \B(p) \ar[rr]^{\psi} \ar[ld]_{\tau_B} \ar'[d][dd]^{\chi} & 
       & S(p) \ar[ld]_{\tau_S} \ar[dd]^{\pi} \\
\B(p) \ar[rr]^(0.7){\psi} \ar[dd]_{\chi} & 
       & S(p) \ar[dd]^(.3){\pi} &\\
& \A(p) \ar[ld]_{\tau_A} \ar'[r]^{\phi}[rr] & 
       & U(p) \ar[ld]_{\tau_U} \\
\A(p) \ar[rr]^{\phi} & & U(p) &
}$$
(recall that $\g=\gl_N$, $S(p)={[S(\g)^{\otimes p}]}^\g$, 
$U(p)={[U(\g)^{\otimes p}]}^\g$; the maps were all defined in the text
above, in particular, $\pi$ is the Poincar\'e-Birkhoff-Witt isomorphism).
\end{lemma}

\begin{proof}
What we are really interested in is the commutativity of the top face
of the cube. Since all $\chi$'s, $\pi$'s and $\tau$'s are linear isomorphisms,
this follows from the commutativity of the remaining 5 faces
which was proved in various places of this paper above.
\end{proof}

The commutativity of the upper face of the cube implies that, to prove the
non-$\tau_B$-invariance of a Jacobi diagram (an element of $\B(p)$) 
it is enough to prove the non-$\tau_S$-invariance of its image under
$\psi$ in the necklace algebra. Since the minimal degree of a 
non-invertible necklace is 6 and we need to find a non-invertible
diagram of an odd degree, the smallest example should be searched 
in degree 7. And, indeed, it exists --- thus providing a second proof of the
Theorem.

\begin{proposition}
The following diagram is non-zero as an element of the space $B(2)$:
\begin{center}
H\quad=\quad\raisebox{-14mm}{\includegraphics[width=3cm]{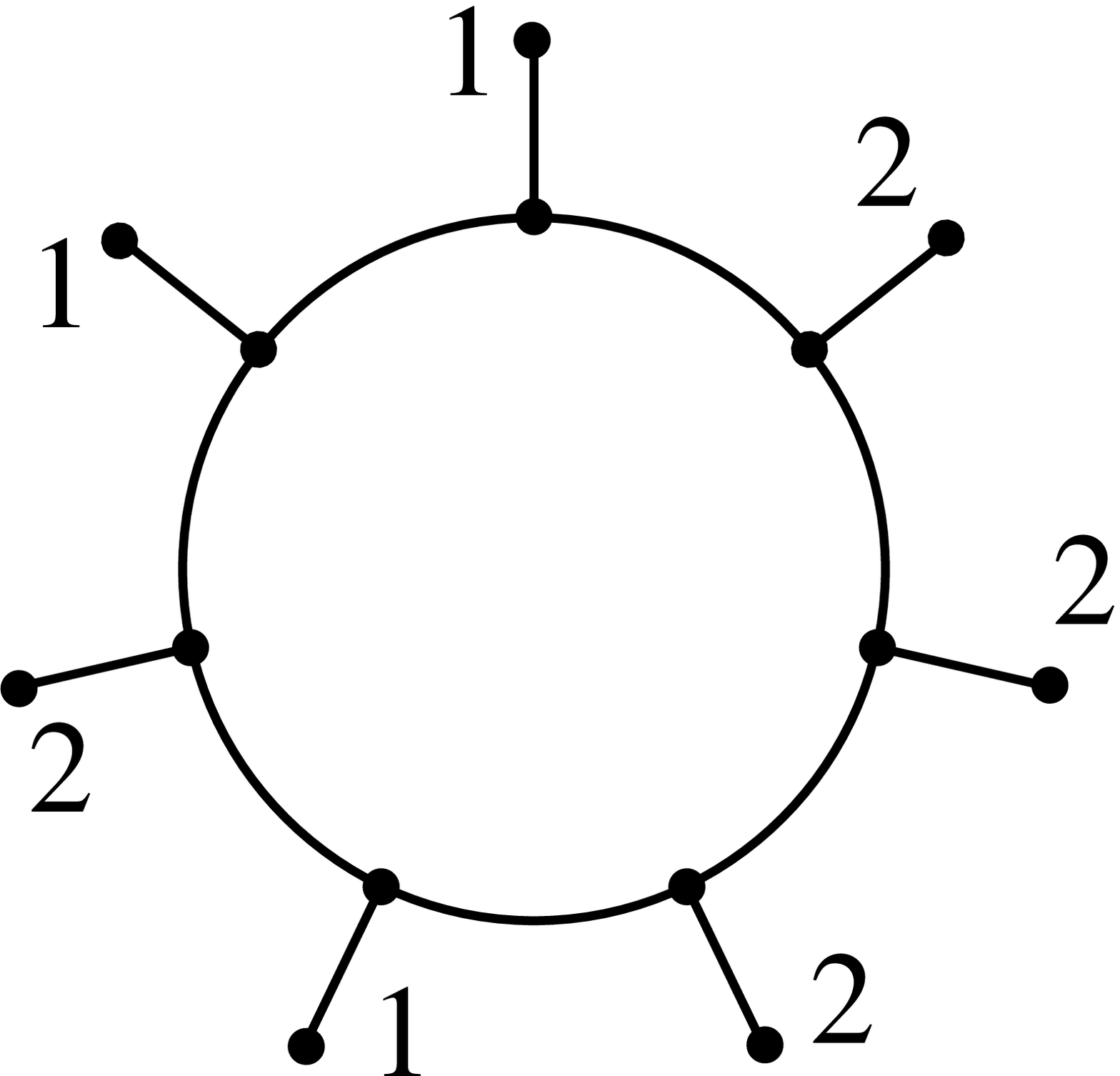}}
\end{center}
\label{prop2}
\end{proposition}

\begin{proof}
A direct computation shows that the image of this diagram in the 
algebra $S(2)$ (for $N$ big enough), when expressed through necklaces,
equals
$$
  N(x_{1121222}-x_{1122212})+3x_2(x_{112212}-x_{112122}),
$$
which is different from zero. The entire expression for $\psi(H)$ involves
128 terms, out of which only 8 contain nonsymmetric necklaces, e.g.
$$
\raisebox{-15mm}{\includegraphics[height=30mm]{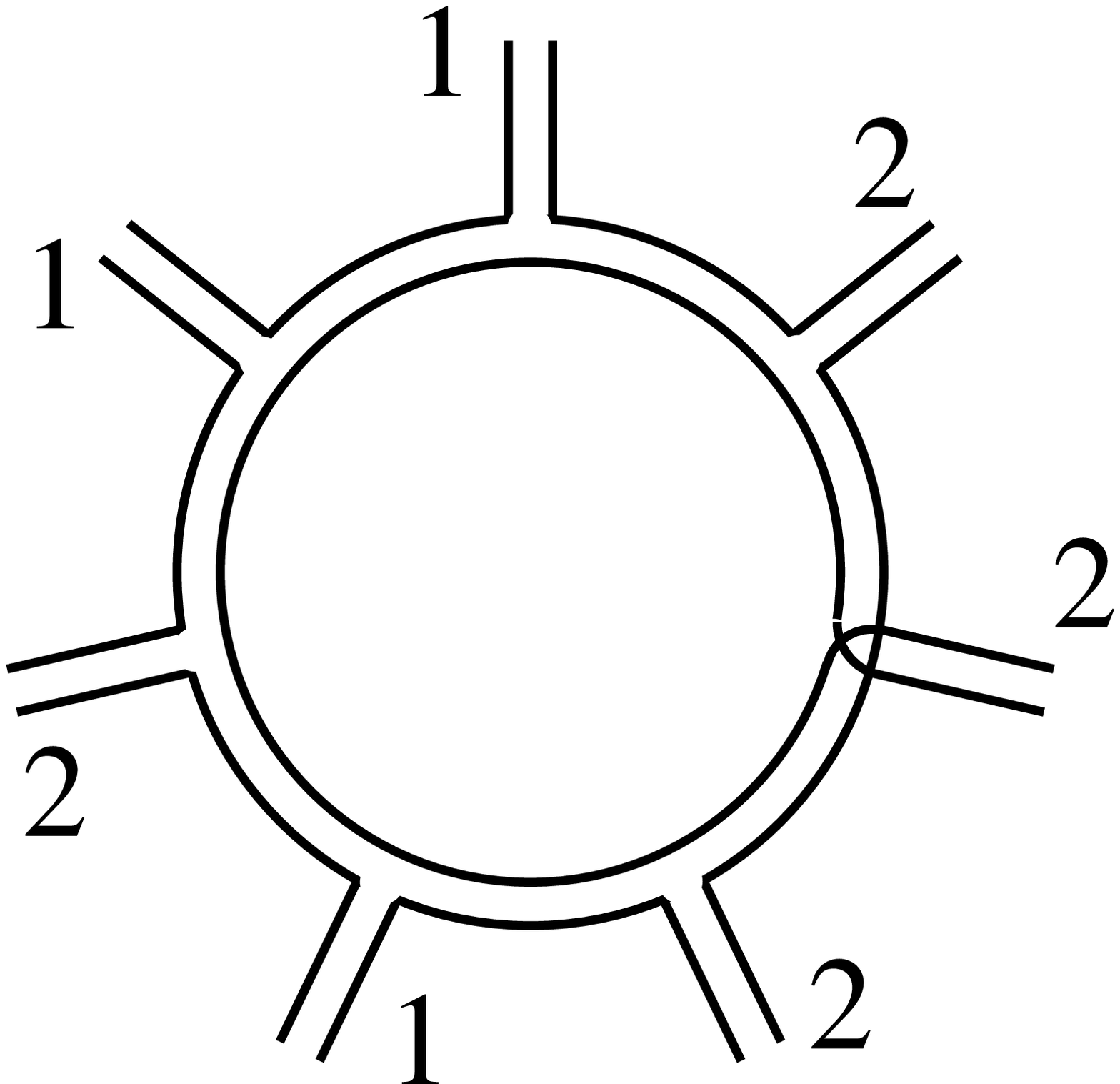}}
\quad \longmapsto \quad -x_2x_{112122}.
$$
The same result can also be obtained by a computer program available at
\cite{Calc}.
\end{proof}

\section{Deframing}
\label{defram}

As we mentioned before, the proofs given in Sections \ref{secphi} and 
\ref{secpsi}, actually refer to the case of \textit{framed} links,
because the 1-term relations were not taken into consideration.
In principle, detecting the orientation of a framed link is easier, because
it contains an additional structure which may not be preserved by the
inversion. Note that an unframed link can be considered as a framed one,
so the Theorem in the unframed case is in any case stronger. Here we explain
why does it hold.

Indeed, let $\A'(2)$ be the space $\A(2)$ modded out by 1-term
relations, i.e. by the ideal generated by the diagrams $a_1$ (with 1 chord on
the first component of the support) and $a_2$ (with 1 chord on the second
component of the support): $\A'(2)=\A(2)/\langle a_1,a_2 \rangle$,
where angular brackets denote the 2-sided ideal with given generators.
The quotient algebra $\A'(2)$ can also be considered as a subalgebra of 
$\A(2)$:
by the structure theorem for cocommutative Hopf algebras, $\A(2)$ is the
universal enveloping algebra over its primitive subspace $P$, so if we take
the subspace of $P$ spanned by all connected diagrams except for $a_1$ and
$a_2$, we will obtain an inclusion $\A'(2)\subset\A(2)$.
Since the diagrams displayed in Proposition \ref{prop1} belong to $\A'(2)$,
this means that non-invertibility holds for unframed long links, too.
\medskip

Now let $\chi_2:\A(2)\to\B(2)$ be the isomorphism discussed above.

\begin{lemma}
The subspace $\chi^{-1}(\A'(2))=\B'(2)$ coincides with the subalgebra of 
$\B(2)$ generated
by all connected Jacobi diagrams except 
$b_1=\raisebox{-0.5mm}{\includegraphics[height=4mm]{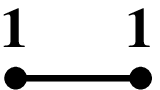}}$ and
$b_2=\raisebox{-0.5mm}{\includegraphics[height=4mm]{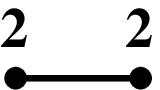}}$.
\end{lemma}

\begin{proof}
The assertion is non-trivial, because $\chi$ does not preserve the
multiplication, and the subalgebras $\A'(2)$ and $\B'(2)$ are described by
their generators in the sense of different multiplications.
However, a little thinking about the action of the operator $\chi$ shows that
its value on a product of Jacobi diagrams different from $b_1$ and $b_2$
is equal to a linear combination of products of connected chord
diagrams different from $a_1$ and $a_2$.
\end{proof}

The lemma shows that the `heptapus' element displayed in Proposition
\ref{prop2} belongs indeed to the subalgebra $\B'(2)$ responsible for
the Vassiliev invariants of unframed 2-component long links. Since it is
non-zero, the existence of a 7 order invariant that can distinguish the
orientation follows.

\section{Related open problems}

\textbf{Problem 1.} Is there a $\tau$-equivariant linear mapping 
$\B(1)\to \B(2)$ whose image 
is not a priori contained in the $\tau$-invariant subspace of $\B(2)$?
\smallskip

\textsl{Remark.} If such a mapping existed, then one could try to apply
the techniques of this paper to the invertibility problem for \textit{knots}.
The standard doubling operator $\Delta:\B(1)\to\B(2)$ (as in \cite{BN1,LM})
is no good, as one can easily see.
\medskip

\textbf{Problem 2.} Do Vassiliev invariants detect invertibility of closed
2-component links?
\smallskip

\textsl{Remark.} The combinatorial object responsible for finite type
invariants of closed links is the space of diagrams on 2 circles, call it
$\A(S^1_2)$. There is an evident epimorphism $\A(\R^1_2)\to\A(S^1_2)$,
but its kernel is non-trivial, so it cannot be used to define a Kontsevich 
type weight system with values in $U^{\otimes2}$.
\medskip

\textbf{Problem 3.} Is 7 the minimal degree of Vassiliev invariants that can
be used to distinguish the orientation of long links with 2 components?
\smallskip

\textsl{Remark.}
The question is reduced to the existence of non-invertible diagrams of
degree smaller than 7.
For example, we do not know if the following diagram 
$$
\raisebox{-12mm}{\includegraphics[height=25mm]{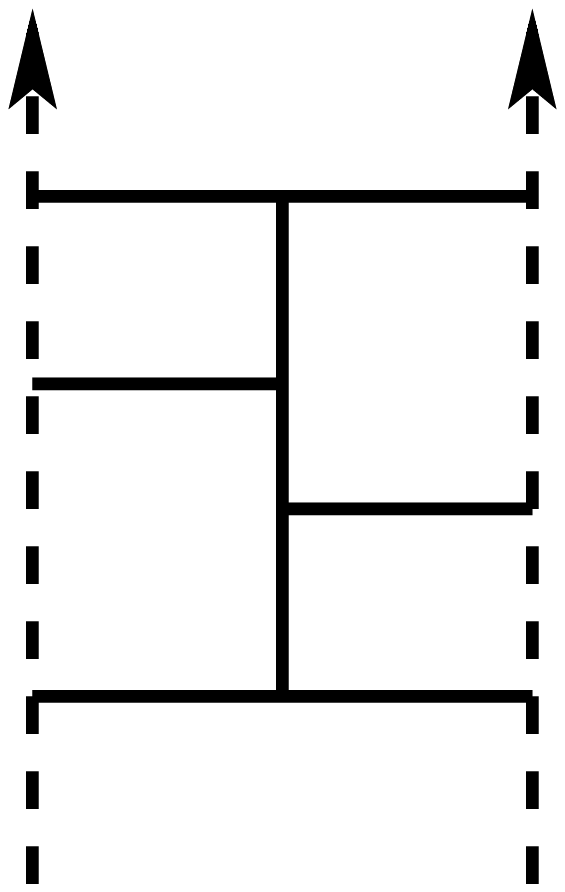}}
\quad-\quad
\raisebox{-12mm}{\includegraphics[height=25mm]{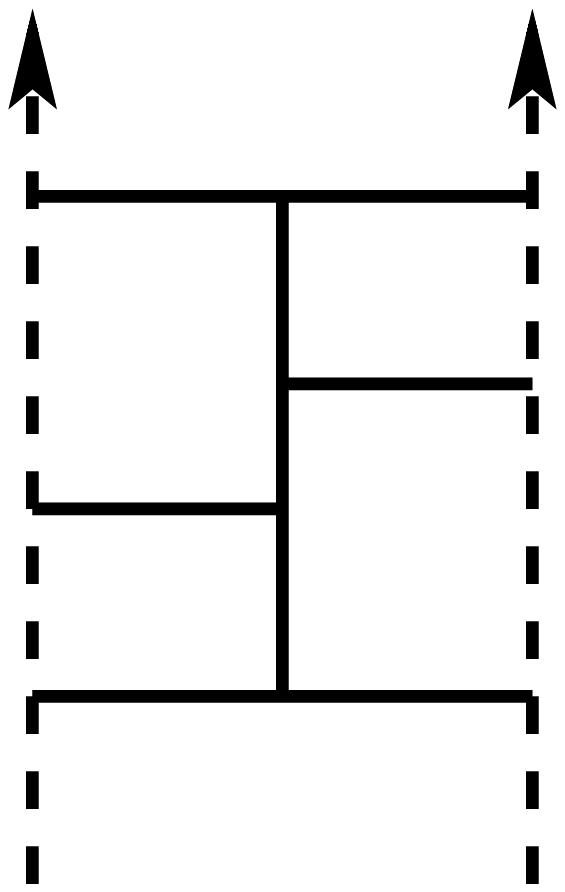}}
\quad=\quad
\raisebox{-12mm}{\includegraphics[height=25mm]{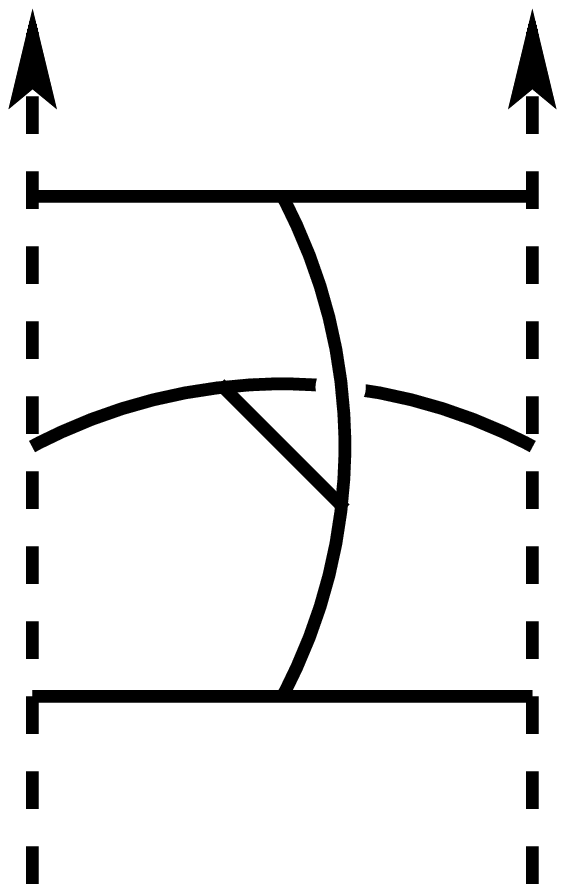}}
$$
is zero or not as an element of $\A_5(2)$.

\section*{Acknowledgements}

We thank A.\,Khoroshkin, S.\,Loktev and E.\,B.\,Vinberg for useful comments
on invariant theory, and S.\,Lando who read the text and made some useful
remarks. The first author was supported by RFBR grant 05-01-01012a.

\rightline{\tt duzhin@pdmi.ras.ru}

\begin{thebibliography}{BN2}

\bibitem[BN1]{BN1} D.\,Bar-Natan, {\it On the Vassiliev knot invariants},
     Topology, {\bf 34} (1995) 423--472.

\bibitem[BN2]{BN2} D.\,Bar-Natan, {\it Vassiliev Homotopy String Link 
    Invariants}, Journal of Knot Theory and its Ramifications 4-1 
    (1995) 13--32.

\bibitem[Calc]{Calc} S.\,Duzhin. \textit{Computer programs and data files 
for the calculation of the weight systems $\phi$ and $\psi$}.
Online at \verb#http://www.pdmi.ras.ru/~arnsem/dataprog/OrLinks/#.

\bibitem[CD]{CD} S.\,V.\,Chmutov and S.\,V.\,Duzhin, \textit{CDBook. 
Introduction to
Vassiliev Knot invariants}, preliminary draft version of a book, 
Online at \verb#http://www.pdmi.ras.ru/~duzhin/papers/#.

\bibitem[Ka]{Ka} A.\,Kawauchi, \textit{The invertibility problem on  
amphicheiral excellent knots}.
Proc. Japan Acad. Ser. A Math. Sci. 55 (1979), no. 10, 399--402.

\bibitem[Kon]{Kon} M.\,Kontsevich, {\it Vassiliev's knot invariants},
    Adv. in Soviet Math., {\bf 16} Part 2 (1993) 137--150.

\bibitem[L1]{Lin1} X.-S.\,Lin,
{\it Finite type link invariants and the invertibility of links}, 
Math. Res. Letters 3 (1996), no. 3, pp. 405--417. 
Online at \texttt{arXiv:q-alg/9601019}.

\bibitem[L2]{Lin2} X.-S.\,Lin,
{\it Finite type link-homotopy invariants}, 
l'Enseignement Math\'ematique 47 (2001), pp. 315--327. 
Online at \texttt{arXiv:math.GT/0012095}.

\bibitem[LM]{LM} T.\,Q.\,T.\,Le, J.\,Murakami, {\it The universal
    Vassiliev--Kontsevich invariant for framed
    oriented links}, Compositio Math. {\bf 102} (1996) 41--64.

\bibitem[Sta]{Sta} T.\,Stanford, {\it Finite-type invariants of knots, links and
graphs}. Topology, Vol. 35, no. 4, pp. 1027--1050.

\bibitem[Tro]{Tro} H.\,F.\,Trotter, {\it Non-invertible knots exist}.
Topology 2 (1964), 275--280.

\bibitem[VP]{VP} 
E.\,B.\,Vinberg and V.\,L.\,Popov, {\it Invariant theory}. In: Itogi Nauki Tekh., 
Ser. Sovrem. Probl. Mat., Fundam. Napravleniya 55 (1989), 137-314. English
transl.: Encyclopaedia Math. Sci. vol 55, Springer, 1994.
% example 1 in no. 9.5, p. 291

\bibitem[Vog]{Vog} P.\,Vogel, {\it Algebraic structures on modules
     of diagrams}, Institut de Math\'{e}\-matiques de Jussieu,  
     Pr\'{e}publication 32, August 1995, Revised in 1997,
 \verb#http://www.math.jussieu.fr/~vogel/#.

\bibitem[W]{Weyl} H.\,Weyl, {\it The Classical Groups: Their Invariants and 
Representations}. Princeton University Press, 1997.

\end{thebibliography}
\end{document}